\documentstyle[12pt]{amsart}
\renewcommand{\marginpar}[1]{}
\catcode`\@=12

\def\Empty{}
\newcommand\oplabel[1]{
  \def\OpArg{#1} \ifx \OpArg\Empty {} \else
  	\label{#1}
  \fi}
		
\long\def\realfig#1#2#3#4{
\begin{figure}[htbp]
\centerline{\psfig{figure=#2,width=#4}}
\caption[#1]{#3}
\oplabel{#1}
\end{figure}}

\newcommand{\comm}[1]{}
    
\input{psfig}
\newtheorem{thm}{Theorem}[section]
\newtheorem{cor}[thm]{Corollary}
\newtheorem{lem}[thm]{Lemma}

\newtheorem{schw}{Schwarz Lemma}

\theoremstyle{definition}

\newcommand{\QED}{\rlap{$\sqcup$}$\sqcap$\smallskip}

\theoremstyle{remark}
\newtheorem{rem}{Remark}[section]

\newcommand{\diam}{\operatorname{diam}}
\newcommand{\dist}{\operatorname{dist}}
\newcommand{\cl}{\operatorname{cl}}
\newcommand{\tl}{\tilde}
\newcommand{\eps}{\epsilon}

\newcommand{\EE}{\cal E}
\renewcommand{\frak}[1]{\EE}

\numberwithin{equation}{section}
\newsymbol\Subset 1362

\newcommand{\secref}[1]{\S\ref{#1}}
\newcommand{\lemref}[1]{Lemma~\ref{#1}}
\newcommand{\corref}[1]{Corollary~\ref{#1}} 
\newcommand{\figref}[1]{Fig.~\ref{#1}}
\newcommand{\ang}[2]{\widehat{(#1,#2)}}
\newcommand{\C}[1]{\Bbb C_{#1}}

\begin{document}

\title[Dynamics of quadratic polynomials: Complex Bounds]%
{Dynamics of quadratic polynomials:\\Complex bounds for real maps }
\author {Mikhail Lyubich }
\address{
Mikhail Lyubich\\
Mathematics Department\\
SUNY at Stony Brook\\
Stony Brook, NY 11794}
\email{
mlyubich@math.sunysb.edu}
\author { Michael Yampolsky}
\address{
Michael Yampolsky\\
Mathematics Department\\
SUNY at Stony Brook\\
Stony Brook, NY 11794}
\email {yampol@math.sunysb.edu}
\thanks{This work was supported in part by Sloan Research Fellowship
and by NSF grants DMS-8920768 and DMS-9022140 (at MSRI).}
\date{\today}

\maketitle

\begin{abstract}
We extend Sullivan's complex a priori bounds to real quadratic
polynomials with essentially bounded combinatorics. Combined with the
previous results of the first author, this yields complex bounds for
all real quadratics. Local connectivity of the corresponding Julia
sets follows.
\end{abstract}

 
\section{Introduction}

Complex a priori bounds proved to be a key issue of the Renormalization 
Theory. They lead to rigidity results, local connectivity of Julia sets
and the Mandelbrot set, and convergence of the renormalized maps (see
\cite{HJ,L,McM,MS,R,S}).

By definition, this property means that the renormalized maps $R^n f$
have fundamental annuli with a definite modulus. For real infinitely
renormalizable maps with bounded combinatorics this property was 
proven by Sullivan ( see \cite{S} and \cite{MS} ). In \cite{L}                                                           %
 complex bounds were
proven for real quadratics of ``essentially big type". The gap in between
\cite{S} 
and \cite{L} consists of maps with ``essentially bounded type".
Loosely speaking this means that a big period of renormalized maps
is created only 
by saddle-node behavior of the return maps. The goal of this paper is to
analyze this specific phenomenon. 

\begin{thm}
\label{bounds}
  Real infinitely renormalizable quadratics with 
essentially bounded combinatorics 
  have complex a priori bounds. 
\end{thm}

This fills the above mentioned gap:

\begin{cor}
\label{bounds1}
\footnote {
 Levin and van Strien have recently announced  
a different proof of this result \cite{LS}.}
 All infinitely renormalizable real quadratics have 
  complex a priori bounds.
\end{cor}

Let us mention here only one consequence of this result.
By  the result of Hu and Jiang \cite{HJ,J}, complex a priori bounds
and one extra combinatorial assumption (see \cite{McM2}) imply
local connectivity of the Julia set $J(f)$. 
On the other hand, the Yoccoz Theorem gives local connectivity
of $J(f)$ for at most finitely renormalizable quadratic maps 
(see \cite{H}, \cite{L1} or \cite{M}). Thus we have

\begin{cor}
\label{lc}
 The Julia set of any real quadratic map is locally
  connected. 
\end{cor}

Theorem \ref{bounds1} is closer to  \cite{S} rather than \cite{L}.
It turns out, however, that Sullivan's Sector Lemma (see \cite{MS}) 
is not valid
 for essentially bounded  (but unbounded) combinatorics: 
The pullback of the plane with two slits is not necessarily contained in a
definite sector. 
\marginpar{Would be good to explain this in the Appendix (when everything else
 is done)} What turns out to be true instead is that the 
{\it little Julia sets} $J(R^n f)$  are contained in
a definite sector. 



We derive this version of the Sector Lemma
 from the following quadratic estimate for the renormalized maps:
\begin{equation} \label{key estimate}
 |R^n f(z)|\geq c|z|^2, 
\end{equation}
with an absolute $c>0$.  
The proof of (\ref{key estimate})
is the main technical concern of this work.  
(By the way, this estimate immediately implies
that the little Julia sets $J(R^n f)$ are commensurable with the
corresponding periodic intervals, which already yields local connectivity
of $J(f)$ at the critical point.)

Let $\mod (f)$ denote the supremum of the  moduli of the fundamental annuli of $f$. 
The work \cite{L}  gives a   
 criterion when   $\mod (R f)$ is big. Let us call the combinatorial
parameter responsible for this
the {\it essential period} $p_e(f)$. Loosely speaking this is the period
of the corresponding periodic interval of $f$ modulo the saddle-node cascades
(see \secref{nest} for the precise definition).

\begin{cor}\label{big space} There is an absolute constant $\gamma>0$ and two 
functions $\mu(p)>\nu(p)>\gamma>0$ going to $\infty$ as $p\to\infty$
  with the following property. Let $f$ be an
infinitely renormalizable quadratic polynomial and $p_n=p_e(R^n f)$. Then
$$\nu(p_n) \leq \mod (R^n f)\leq \mu(p_n). $$ 
\end{cor}

Let us briefly outline  the structure of the paper. 
\secref{prelim} contains some background and technical preliminaries.  
In \secref{outline} we state the main technical lemmas, and derive from them
our results.
In \secref{bounded}
 we give a quite simple proof
of complex bounds in the case of bounded combinatorics, which will model
the following argument. In \secref{nest} essentially
bounded combinatorics is described. In the next section, \secref{cascades},
 saddle-node cascades are analyzed.  The final section, \secref{inductive},
 contains the proof of the 
main technical lemmas. 

\medskip{\it Remark 1.} Theorem 1.1 allows a straightforward extension onto higher
degree unimodal polynomials. 

{\it 2.} This paper is a part of series of notes on dynamics of quadratic polynomials,
  see \cite{L4}. 

\medskip{\bf Acknowledgment.} The authors thank MSRI where part of this work
  was done for its hospitality.

\section{Preliminaries}
\label{prelim}
\subsection{General notations and terminology} 
\marginpar{Put here all sort of stuff}
We use $|J|$ for the length of an interval $J$, 
 $\dist$ and $\diam$ for the Euclidean distance and diameter in $\Bbb C$.
Notation $[a,b]$ stands for  the (closed) interval with endpoints $a$ and $b$ without
specifying their order. 

Two sets  $X$ in $Y$ in $\Bbb C$ are called {\it $K$-commensurable}
or simply {\it commensurable} if 
$$K^{-1}\leq \diam X/\diam Y\leq K$$ with a constant $K>0$ which 
may depend only on the specified combinatorial bounds.

We say that an annulus $A$ has a {\it definite modulus}
 if $\mod A\geq \delta>0$,
where $\delta $ may also depend only on the specified combinatorial bounds.  

For a pair of intervals $I\subset J$ we say that $I$ is contained 
{\it well inside} of $J$ if for any of the components $L\subset J\setminus I$,
 $|L|\geq K|I|$ where the constant $K>0$ may depend only on the specified
quantifiers.

A smooth interval  map $f: I\rightarrow I$ is called {\it unimodal} if it has a single
critical point, and this point is an extremum. 
A $C^3$ unimodal map is called {\it quasi-quadratic} if it has negative Schwarzian derivative,
and  its critical point is non-degenerate. 

Given a unimodal map $f$ and a point $x\in I$, $x'$ will denote the dynamically
symmetric point, that is, such that $fx'=fx$. Notation $\omega(z)$ means as usual the
limit set of the forward orbit $\{f^n z\}_{n=0}^\infty$. 

\subsection{Hyperbolic disks}
Given an interval $ J\subset \Bbb R$, let
$\C{J}\equiv \Bbb C\backslash (\Bbb R\backslash J)$ denote
the plane slit along two rays. Let  ${\bar {\Bbb C}}_J$
 denote the completion of
this domain in the path metric in $\C{J}$ (which means that we add to $\C{J}$
 the banks of the slits).

By symmetry, $J$ is a hyperbolic geodesic in $\C{J}$.
The {\it geodesic neighborhood} of $J$  of radius
$r$ is the 
set of all points in $\C{J}$
 whose hyperbolic distance to $J$ 
 is less than $r$. It is easy to see that such a neighborhood is
the union of two $\Bbb R$-symmetric segments of Euclidean disks 
based on $J$
 and having angle $\theta=\theta(r)$ with $\Bbb R$. Such a hyperbolic 
disk will be denoted by $D_{\theta}(J)$ (see Figure 1).
Note, in particular, that the Euclidean disk $D(J)\equiv D_{\pi/2}(J)$
can also be interpreted as a hyperbolic disk.
%

These hyperbolic neighborhoods were introduced into the subject by Sullivan
[S]. They are a key tool for getting complex bounds due to
the following version of the Schwarz Lemma:

\begin{schw}
  Let us consider two intervals  
$J'\subset J\subset {\Bbb R}$. 
Let $\phi:\C{J}\rightarrow\C{J'}$ be an analytic map such that
$\phi(J)\subset J'$. Then for any $\theta\in (0, \pi),$ 
$\phi (D_\theta (J))\subset D_\theta (J')$.
\end{schw}

Let $J=[a,b]$.  
For a point $z\in {\bar {\Bbb C}}_J$, the 
{\it angle between $z$ and 
 $J$}, $\ang{z}{J}$
is the least of the angles  between the intervals $[a,z]$, $[ b, z]$  and 
the corresponding rays $(a,-\infty]$, $[b, +\infty)$ of the real line,
measured in the range $ 0\leq \theta\leq \pi$. 

 We will use the following observation to control the expansion of the inverse
branches. 

\begin{lem}\label{goodangle}
Under the circumstances of the Schwarz Lemma, let us consider 
 a point $z\in \C{J}$ such that $\dist(z, J)\geq |J|$ and 
$\ang{z}{J}\geq\epsilon$. Then
$${\dist(\phi z, J')\over |J'|}\leq C{\dist(z, J)\over |J|}$$ 
for some constant $C=C(\epsilon)$
\end{lem}

\realfig{gdang}{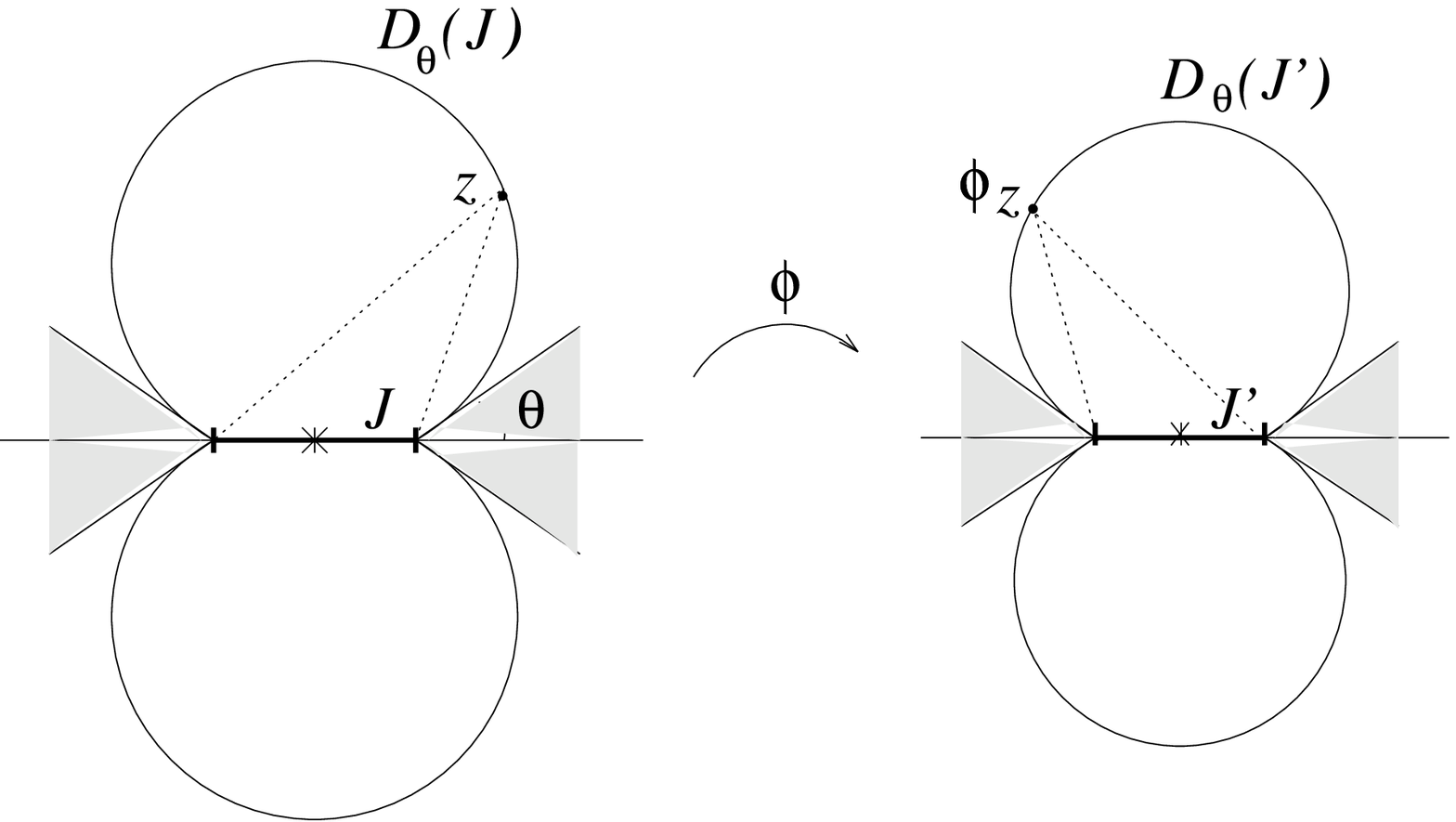}{}{15cm}

\begin{pf} Let us normalize the situation in this way: $J=J'=[0,1]$.
Notice that the smallest (closed) geodesic neighborhood 
$\cl D_{\theta}(J)$ enclosing $z$ satisfies:\\ 
$\diam D_{\theta}(J)\leq C(\epsilon)\dist (z, J)$ (cf \figref{gdang} ).

Indeed, if $\theta\geq \epsilon/2$ then 
$\diam D_{\theta} (J)\leq C(\epsilon)$, which is fine since
$\dist (z, J)\geq 1$.\\
Otherwise the intervals $[0,z]$ and $[1,z]$ cut out sectors of angle 
size at least $\epsilon/2$ on the circle $\partial D_{\theta}(J)$.
Hence the lengths of these intervals are commensurable with 
$\diam D_{\theta}(J)$ (with a constant depending on $\epsilon$).
Also, by elementary trigonometry these lengths are at least 
$\sqrt{2}\dist(z, J)$, provided that $\dist (z, J)\geq |J|$.

By Schwarz Lemma, 
 $\dist(\phi z, J')\leq\diam( D_{\theta}(J')),$ and the claim follows.
\end{pf}

\subsection{Square root}
In the next lemma we collect for future reference some elementary properties of the
square root map. Let $\phi(z)=\sqrt{z}$ be the  branch of the square
root mapping the slit plane ${\Bbb C}\setminus {\Bbb R}_-$ into itself.

\begin{lem}
\label{square root}
 Let $K>1$, $\delta>0$, 
$K^{-1}\leq a\leq K$, $T=[-a, 1]$, $T'=[0,1]$. Then:
\itemize{

\item $\phi D_\theta(T)\subset D_{\theta'}(T')$, 
with $\theta'$ depending on $\theta$ and $ K$ only.

\item 
 If $z'\in \phi D(T)\setminus D([-\delta,1+\delta])$, then 
$$\ang{z'}{T'}>\epsilon(K, \delta)>0\quad and\quad
C(K,\delta)^{-1}<\dist({z'},{T'})< C(K, \delta).$$}
\end{lem}

\begin{lem}
\label{square root2}
 Let $\zeta\in \Bbb C$,  $J=[a,b]\subset [0,+\infty)$.
 $\zeta'=\phi(\zeta)$, $J'=[a',b']=\phi J$.  Then: \itemize{

\item If $\dist(\zeta,J)>\delta |J|$ then
$$\frac{\dist(J',\zeta')}{|J'|}<C(\delta)\frac{\dist(J,\zeta)}{|J|}.$$

\item  Let $\theta$ denote the angle between $[\zeta, a]$ and the
 ray of the real line which does not contain $J$; $\eta'$ denote the angle
between $[\zeta', b']$ and the corresponding ray of the real line.  
If $\theta\leq \pi/2$ then $\eta'\geq \pi/4$.}  
\end{lem}

(According to our convention, in the last statement we don't assume that
$a<b$.)

\subsection{Branched coverings} 
Let $0\in U'\subset U\subset {\Bbb C}$ 
be two topological disks different from the whole plane,
and $f: U'\rightarrow U$ be an analytic branched double covering map with critical point at 0. Thinking of it as a dynamical system, one can naturally
define the {\it filled Julia set} $K(f)$ and {\it the Julia set} $J(f)$.
 Namely, the filled Julia set is the set of non-escaping points,
$$K(f)=\{z: f^n z\in U', n=0,1,\ldots\},$$ and $J(f)=\partial K(f)$.
These sets are not necessarily compact.

If additionally $\cl U'\subset U$ then the map $f$ is called 
{\it quadratic-like}. The Julia set of a quadratic-like map is compact,
and  this is actually the criterion
for being quadratic-like (for appropriate choice of
domains):

\begin{lem} [compare \cite{McM2}, Proposition 4.10]
 \label{pol-like} Let $U'\subset U$ be two topological disks.
 and  $f: U'\rightarrow U$  
be a double branched covering with non-escaping critical 
point and compact Julia set. Then there are topological discs
 $J(f)\subset V'\subset V\subset U$
 such that the restriction $g: V'\rightarrow V$  is 
quadratic-like.  Moreover, if $\mod (U\setminus K(f)\geq \epsilon>0$ then
$\mod (V\setminus V')\geq \delta(\epsilon)>0$.
\end{lem}

\begin{pf} Let us consider the topological annulus
 $A=U\setminus K(g)$. Let $\phi: A\rightarrow R=\{z: 1<|z|<r\}$ be
its uniformization by a round annulus. It conjugates $g$ to a map
$G: R'\rightarrow R$ where $R'$ is a subannulus of $R$ with the same
inner boundary, unit circle $S^1$. As $G$ is proper near the unit circle,
it is continuously extended to it, and then can be reflected to the
symmetric annulus. We obtain the double covering 
map $\hat G: \hat R'\rightarrow \hat R$ of the symmetric annuli
preserving the circle. Moreover $\hat R$ is a round annulus of modulus
at least $2\epsilon$.

Let $l$ denote the hyperbolic length on $\hat R$, $\hat V$ denote
the hyperbolic 1-neighborhood of $S^1$, and $\hat V'=\hat G^{-1} \hat V\subset
\hat V$.
 As $\hat G: S^1\rightarrow S^1$ is a double covering, we have:
$$2 l(S^1)= \int_{S^1} \|Df(z)\| dl\leq \max_{S^1}\|Df(z)\| l(S^1),$$
so that $\max_{S^1} \|Df(z)\|\geq 2$. As $\mod \hat R\geq 2\epsilon$,
$l(S^1)\leq L(\epsilon)$. Hence $\| Df(z)\|\geq \rho(\epsilon)>1$
for all $z\in \hat V$. It follows that $\hat V'$ is contained in 
$(1/\rho(\epsilon))$-neighborhood of $S^1$. But then each component
of $V\setminus V'$ is an annulus of modulus at least $\delta(\epsilon)>0$.

We obtain now the desired domains by
going  back to $U$:  $V=\phi^{-1} \hat V$,
$V'=\phi^{-1} \hat V'.$  
\end{pf}

 Let us supply the space ${\cal B}$ of double branched maps considered above
with the {\it Caratheodory topology} (see \cite{McM}). Convergence 
of a sequence $f_n : U_n'\rightarrow U_n$  in this
topology means Caratheodory convergence
of $(U_n,0)$ and $U_n', 0)$,
 and compact-open convergence of $f_n$.

\subsection{Epstein class} A double branched map $f: U'\rightarrow U$ of class
${\cal B}$ belongs to Epstein class if $U=\C{T}$, $U'$ is an $\Bbb R$-symmetric
domain meeting the  real line along an interval $T\supset T'$, and the map $f$
is $\Bbb R$-symmetric. In this case its restriction $f: T'\rightarrow  T$
is a unimodal map. We always normalize $f$ in such a way that 0 is its
critical point.

Given a $\lambda\in (0,1)$, 
let  $\EE_{\lambda}$ denote the space of maps of Epstein class
with 
$$\lambda |T'|\leq |T|\leq \lambda^{-1} |T'|,$$ modulo affine conjugacy
(that is, rescaling of $T$). 

\begin{lem}\label{compactness} For each $\lambda\in (0,1)$, the space 
  $\EE_{\lambda}$  is compact.
\end{lem}

\begin{pf} Normality argument. 
\end{pf}

{\it All maps in this paper will be assumed to belong to some Epstein class.}

 


\subsection{ Renormalization.}
\label{renorm}
We assume that the reader is familiar with the notion of renormalization
in one-dimensional dynamics (see e.g., [MS]).     

Let $f$ be infinitely renormalizable.
Let $P^k\ni 0$ be the central periodic interval corresponding
to the $k$-fold renormalization $R^k f$ of $f$, $n_k$ be its period:
 $f_k\equiv R^kf\equiv f^{n_k}:P^{k}\rightarrow P^k$.
Set $P_m^k=f^m P^k$. 
 We say that the intervals $P^k_i,\; i=0,1,\ldots, n_k-1,$ 
form the {\it cycle of level  $k$ }. 

Note that the periodic interval
$P^k$ is not canonically defined. The maximal choice is
 $P^k=B^k=[\beta_k, \beta_k']$ where $\beta_k$ is the fixed
point of $f_k$ with positive multiplier. The minimal choice is
$P^k=[f_k 0, f_k^2 0]$. 

Let $p_k=n_k/n_{k-1}$ be relative periods.
 {\it Combinatorics} of $f$ is said to be {\it bounded}
 if the sequence of relative periods is bounded.  
Let $G^k_l$ be the {\it gaps} of level $k$, that is the components of
$P^{k-1_j}\setminus\cup P^k_i$.
{\it Geometry} of  $f$ is said to be {\it bounded}
 if there is a $\Delta>0$ and a choice of periodic intervals $P^k_i$, such that 
for any $P_i^{k}, G^{k}_l\subset P_j^{k-1},\;\;$ 
$|P_i^{k}|/|P^{k-1}_j|>\Delta$
and $|G^{k}_l|/|P^{k-1}_j\geq \Delta$.  In other words, all intervals and gaps
of level $k$ contained in some interval of level $k-1$ are commensurable
with the latter.

\proclaim Theorem A. Infinitely renormalizable 
maps with bounded combinatorics have bounded geometry.

For a proof the reader is referred to \cite{G,BL1,BL,S,MS}.

Let $S^k\supset P_1^k$ be  
the maximal symmetric interval around $0$ such that the restriction of
$f_k$ to it is unimodal, and  $T^k=f_k S^k$. 
Then $P^k\subset\ S^k\subset T^k$, and there is a definite
space in between any two of these intervals.
In the case of bounded (and essentially bounded)
combinatorics  all three intervals are  commensurable.
Moreover, if $f$ belongs to Epstein class, then  the renormalizations $f_k$ 
are also maps of  Epstein class, with range $\C{T^k}$.

\proclaim Corollary B. If $f$ is an infinitely renormalizable map
of Epstein class with bounded combinatorics, then all renormalizations 
$R^n f$ belong to some Epstein class   ${\frak{E}}_{\lambda}$.
Hence the sequence $R^n f$ is pre-compact.



\section{Outline of the proof}\label{outline}

\subsection{Main lemmas} 
Let $P^k$, $f_k\equiv R^k f\equiv f^{n_k}$ be as above.

Let us consider the decomposition:
\begin{equation}\label{decompose}
f_k=\psi_k\circ f,
\end{equation}
where $\psi_k$ is  a univalent map from a neighborhood of $P^k_1$ onto $\C{T^k}$.

At \secref{nest} we will define the essential period $p_e(f)$. For the
time being the reader can just replace this by the period
$p(f)=n_1$.

\begin{lem}\label{contraction}
Let $f$ be a $k$ times  renormalizable quadratic map.
Assume that $p_e(R^l f)\leq \bar p$ for $l=0,1,\ldots,k-1$.  
Then there exist constants $C, D$, depending on $\bar p$ only, such that
$\forall z\in \C{T^k}$ with $\dist(z, P^k)\geq |P^k|$ the following 
estimate holds:
\begin{equation}
\label{linear}
 \frac{\dist(\psi_k^{-1} z, P^k_1)}{|P^k_1|}\leq 
C\left(\frac{\dist (z, P^k)}{|P^k|}\right)  +D,
\end{equation}
where $\psi_k$ is the univalent map from \ref{decompose}
\end{lem}


Thus the maps $\psi_k^{-1}$ have at most linear growth depending only on the
combinatorial bound $\bar p$.  
 
Note that if $\ang{z}{P^{k}}>\epsilon>0 $, the
inequality \ref{linear} follows directly from \lemref{goodangle}, with the 
constants depending on $\epsilon$.
Our strategy of proving \lemref{contraction} is to monitor the inverse orbit
of a point $z$ together with the interval $P^{k}$ until they satisfy this
"good angle" condition.

\lemref{contraction}  immediately yields the key quadratic estimate \ref{key estimate},
which in turn implies:
\begin{cor}
\label{contr2}
The little Julia set $J(R^k f)$ is commensurable
with the interval $P^k$.
\end{cor}

Carrying the argument for \lemref{contraction} further, we will prove the following result:

\begin{lem}\label{sector} Under the circumstances of the previous lemma, 
the little Julia set $J(R^k f)$ is contained in the hyperbolic disk
$D_\epsilon ( B^{k}),$ where $\epsilon>0$ depends only on $\bar p$.
\end{lem}

\subsection{Proof of the main results.} 

\medskip\noindent{\it Proof of Theorem \ref{bounds}.} It follows immediately 
from Lemma \ref{sector} and Lemma \ref{pol-like}. \QED


\medskip\noindent{\it Proof of Corollary \ref{bounds1}.}
 By \cite{L}, there is a $\bar p$
such that $\mod (Rf)\geq \mu>0$ for all renormalizable maps $f$
 of Epstein class with $p_e(f)\geq \bar p$. 

So given a quadratic polynomial, we have complex bounds for all
renormalizations $R^{n+1} f$ such that $p_e(R^n f)\geq \bar p$.
For all intermediate levels we have bounds by Theorem \ref{bounds}. \QED


By a {\it puzzle piece} we mean a topological disk bounded by rational external rays
and equipotentials (compare \cite{H,L4,M})

\medskip\noindent

\noindent{\it Proof of Corollary \ref{lc}.}  
By \corref{contr2} the little Julia sets $J(f_k)$
shrink to the critical point. By the Douady and Hubbard renormalization construction 
(see \cite{D,L4,M2}), 
each little Julia set is the intersection of
a nest of puzzle pieces. As each of these pieces contains
a connected part of the Julia set, $J(f)$ is locally connected at the 
critical point.

Let us now prove local connectivity at any other point $z\in J(f)$ (by a standard "spreading around"
argument). Take a puzzle piece $V\ni 0$. The set of points which never visit $V$,
 $Y_V=\{\zeta: f^n\zeta\not\in V, \; n=0,1,\dots\}$, is expanding. (Cover this set by finitely many non-critical 
puzzle pieces, thicken them a bit, and use the fact the branches of the inverse map are contracting with respect
to the Poincar\'{e} metric in these pieces). It follows that if $z\in Y_V$ then there is a nest of puzzle pieces
shrinking to $z$, and we are done.

By \lemref{sector}, there is a nest of puzzle pieces $V^k\supset J(f_k)$ contained in the Poincar\'{e} disk
$D_{\theta}(B^k)$, with $\theta>0$ depending only on $\bar p$.  But because of bounded geometry
 (or, more generally, "essentially bounded geometry", see \S 5),   there is a definite gap  between the 
interval $B^k$ and the rest of the postcritical set $\omega(0)$. (That is, there is an $\eps=\eps(\bar p)>0$
such that the interval $(1+\eps)B^k\setminus B^k$ does not intersect $\omega(0)$.) Thus the
annuli $R^k=\C{(1+\eps)B^k}\setminus D_\theta(B^k)$ don't meet the postcritical set. Moreover, 
all these annuli are similar and hence have the same moduli.

Assume now that $f^{l_k}z\in V^k$. Then there exist single-valued
 inverse branches $f^{-l_k}: \C{(1+\eps)B^k}\rightarrow \Bbb C$  whose images contain $ z$. 
By the Koebe theorem, they have a bounded 
distortion on puzzle pieces $V^k$. As $U^k=f^{-l_k} V^k$ cannot contain a disk of a definite radius,
we conclude that  $\diam U^k\to 0$. This is the desired nest of puzzle pieces about $z$. \QED

\medskip\noindent{\it Proof of \corref{big space}.} 
This result follows from Theorem D of \cite{L} and Theorem 1.1.\QED


\section{Bounded Combinatorics}
\label{bounded}
We first show the existence of the complex bounds in the case when
the map $f$ has bounded  combinatorics. The result is well-known in
this case \cite{MS,S}, but we give a quite simple proof which will
be then generalized for the case of essentially bounded combinatorics.

\subsection{The $\epsilon$-jumping points}\label{jumping} 

Given an interval $T\in \Bbb R$ let $f: U'\rightarrow\C{T}$ be a 
map of Epstein class.

For a  point $x\in {\Bbb R}\cap U'$ which is not critical for $f^n$,
 let $V_n(x)\equiv V_n(x, f)$ denote the maximal domain containing $x$
which is univalently mapped by $f^n$ onto $\C{T}$. Its intersection with
the real line is the monotonicity interval $H_n(x)\equiv H_n(x, f)$ of $f^n$ containing $x$.
Let  $f^{-n}_x: \C{T}\rightarrow V_n(x)$ denote the corresponding inverse branch
 of $f^{-n}$ (continuous up to the boundary of the slits, with different
values on the different banks). 
If $J$ is an interval on which $f^n$ is monotone, then
the notations $V_n(J)$ and $H_n(J)$ and $f^{-n}_J$ make an obvious  sense.

Take an $x\in \Bbb R$ and a $z\in \C{T}$. If we have a backward orbit of 
$x\equiv x_0, x_{-1}, \dots, x_{-l} $ of $x$ which does not contain 0, 
the {\it corresponding } backward orbit
$z\equiv z_0, z_{-1}, \dots, z_{-l}$ is obtained by applying the appropriate
branches of the inverse functions: $z_{-n}=f_{x_{-n}} z$. The same
terminology is applied when we have a monotone pullback
 $J\equiv J_0,\dots, J_{-l}$ of an interval $J$.

Let $H\supset J$ be two intervals. Let $S_{\theta,\epsilon}(H, J)$ denote the
union of two $2\epsilon$-wedges with vertices at $\partial J$ (symmetric 
with respect to the real line) cut off by  the neighborhood
$D_\theta (H)$ (cf. \figref{fig12}).
\realfig{fig12}{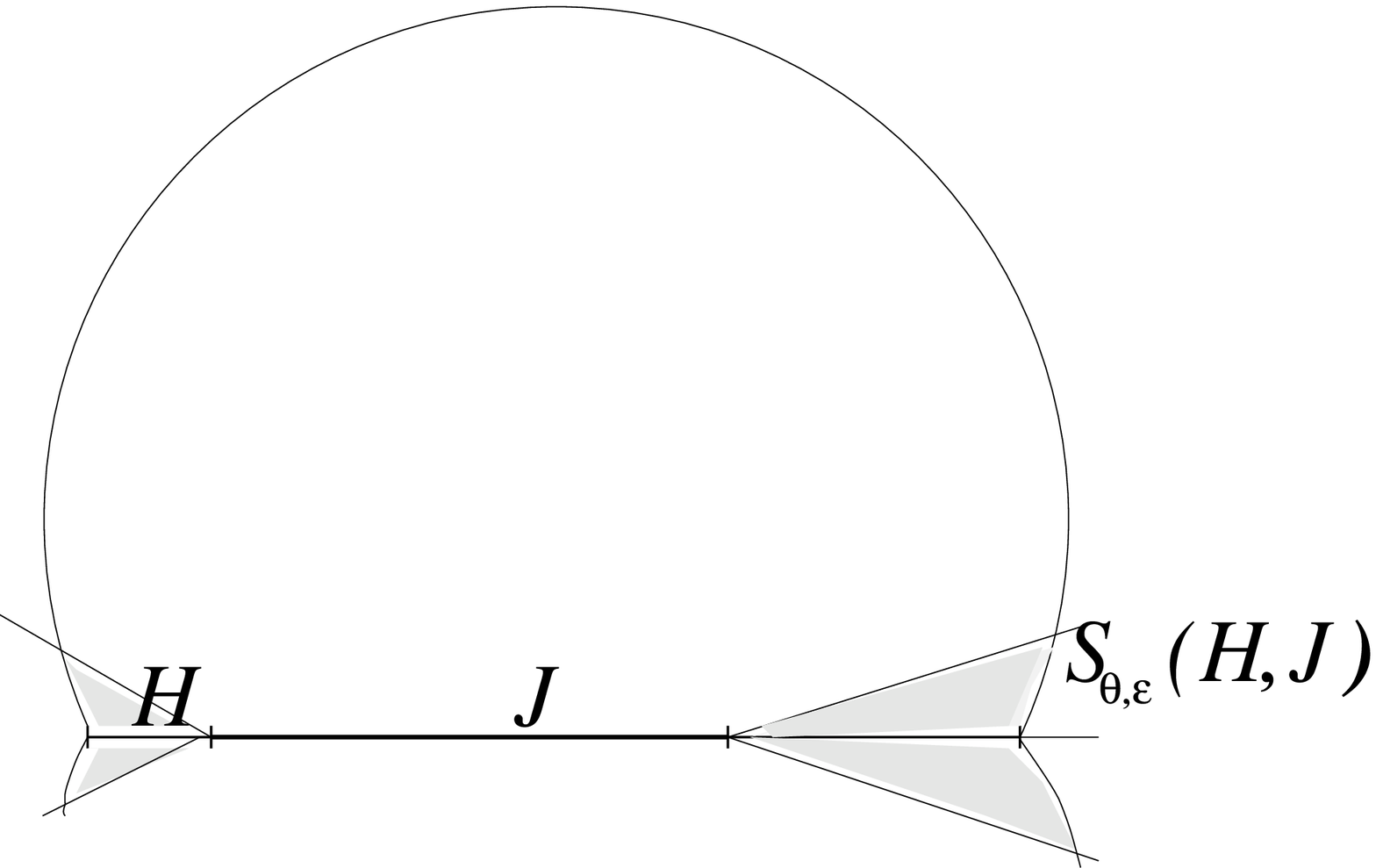}{}{11cm}
Let $Q_{\epsilon}(J)$ denote the complement of the above two wedges
(that is, the set of points looking at $J$ at an angle at least $\epsilon$).

\begin{lem} 
\label{joc}
 Let $f$ be a quadratic map.
 Let $J\equiv J_0, J_{-1},\dots, J_{-l}\equiv  J'$ be a monotone 
pullback of an interval $J$, $z\equiv z_0, z_{-1},\ldots, z_{-l}\equiv z'$ 
be the corresponding backward orbit of a point $z\in \C{T}$. 
Then  for all sufficiently small $\eps>0$ (independent of $f$), 
either $ z_{-k}\in Q_{\epsilon}(J_{-k})$ at some moment $k\leq l$, 
or $z'\in S_{\theta,  \epsilon}(H_l(J'), J')$ with $\theta=\pi/2-0(\epsilon)$.
\end{lem}
If the first possibility of the lemma occurs we say that the backward orbit
of $z$ "$\epsilon$-jumps".
\begin{pf}
Assume that the backward orbit of $z$ does not ``$\epsilon$-jump", that is, $z_{-k}$ belongs
to an $\Bbb R$-symmetric $2\epsilon$-wedge centered at $a_{-k}\in \partial J_{-k}$,
$k=0,1\ldots,l$. By the second statement of \lemref{square root2},
 $f a_{-(k+1)}=a_{-k}$.
Let $M_{-k}=f^{l-k} H_n(J')$, and $b_{-k}$ be the boundary point of $M_k$
on the same side of $J_{-k}$ as $a_{-k}$.
 Let us take the moment $k$ when $b_{-k}=0$. 
At this moment the point $z_{-k}$ belongs to a right triangle based upon
$[a_{-k}, b_{-k}]$ with the $\epsilon $-angle at $a_{-k}$ and the right angle
at $b_{-k}$. Hence
$z_{-k}\in D_{\theta} (M_{-k})$ with $\theta=\pi/2-0(\epsilon)$.
It follows by Schwarz Lemma  that $z'\in D_{\theta} (H_l(J'), J')$, 
and we are done.  
\end{pf}

 Let us state for the further reference in 
\secref{inductive} a straightforward extension 
of the above lemma onto maps of Epstein class:

\begin{lem}\label{joc-2} The conclusion of \lemref{joc} still holds
for all $\eps<\eps(\lambda)$, provided
 $f: U'\rightarrow \C{T}$ is a map of Epstein class $\EE_{\lambda}$,
$U'\subset \C{\lambda T}$, and $z\in D(\lambda T)$.
\end{lem}

\subsection{Proof of \lemref{contraction} (for bounded combinatorics).}

For technical reasons we consider a new family of intervals $\tl{S}^k$
and $\tl{T}^k$, for which $P^k\subset\tilde{S}^k\subset S^k\subset \tl{T}^k
\subset T^k$, each of the intervals is commensurable with the others and
contained well inside the next one, and $f_k (\tl{S}^k)=\tl{T}^k$.

Let us fix a level $k$, and set $n\equiv n_k$, 
\begin{equation}\label{J-orbit}
J_0\equiv P^k_0, J_{-1}\equiv P^k_{n-1},\ldots, J_{-(n-1)}\equiv 
P^k_1.
\end{equation} 
 Take now any point $z_0\in \C{T^k}$ with 
$\dist (z_0, J_0)\geq |J_0|$, and let 
$z_{-1}, \ldots, z_{-(n-1)}$ be its backward orbit corresponding
to the above backward orbit of $J_0$. Our goal is to prove that

\begin{equation}\label{goal}
{\dist(z_{-(n-1)}, J_{-(n-1)})\over |J_{-(n-1)}|}\leq C(\bar p) 
  {\dist (z_0, J_0)\over |J_0|}.
\end{equation}

Take a big quantifier $\bar K>0.$ 
Let is say that $s$ is a "good" moment of time if  $J_{-s}$ is $\bar K$-commensurable with $J_0$.  For example, let $J_{-s}\subset P^l$ and $s\leq n_{l+1}$.
In other words, $s$ is a moment of backward return to $P^l$ preceding the first return to $P^{l+1}$.
By bounded geometry, this moment is good, provided $\bar K$ is selected sufficiently big.

Let us first consider  the initial piece of the orbit, $z_0, \ldots, z_{-n_1},$
 corresponding to the renormalization cycle 
of level $1$. By the first statement of \lemref{square root2},
 for all $s\in [0, n_1-1]$,
\begin{equation}\label{sqrt} 
{\dist(z_{-s}, J_{-s})\over |J_{-s}|}\leq C_0(\bar p) 
  {\dist (z_0, J_0)\over |J_0|}.
\end{equation}

 By \lemref{joc}, either 
\begin{equation}\label{initial}
z_{-n_1}\in D(S^1)\subset D(\tl{T}^1),
\end{equation}
or there is a moment $-s\in [ -n_1, 0]$ when the backward orbit 
$\epsilon$-jumps: $\ang{z_{-s}}{J_{-s}}>\epsilon$. In the latter case
the desired estimate (\ref{goal}) follows from (\ref{sqrt}) and 
\lemref{goodangle}.
In the former case we will proceed inductively:

\begin{lem}
\label{bounded1} Let $J=J_{-s}$ and $J'=J_{-(s+n_l)}$ be two consecutive returns of the
  backward orbit (\ref{J-orbit})  to a periodic interval $P^l$, $l< k$. 
Let $z$ and $z'$ be the corresponding points of the backward orbit of $z_0$.
If $z\in D(\tl{T}^{l})$ then $\dist(z', J')\leq C(\bar p) |\tl T^l|$. Moreover,
either $z'\in D(\tl{T}^{l})$, or  $\ang{z'}{J'}>\epsilon(\bar p)>0$.
\end{lem} 

\begin{pf}

\realfig{fig10}{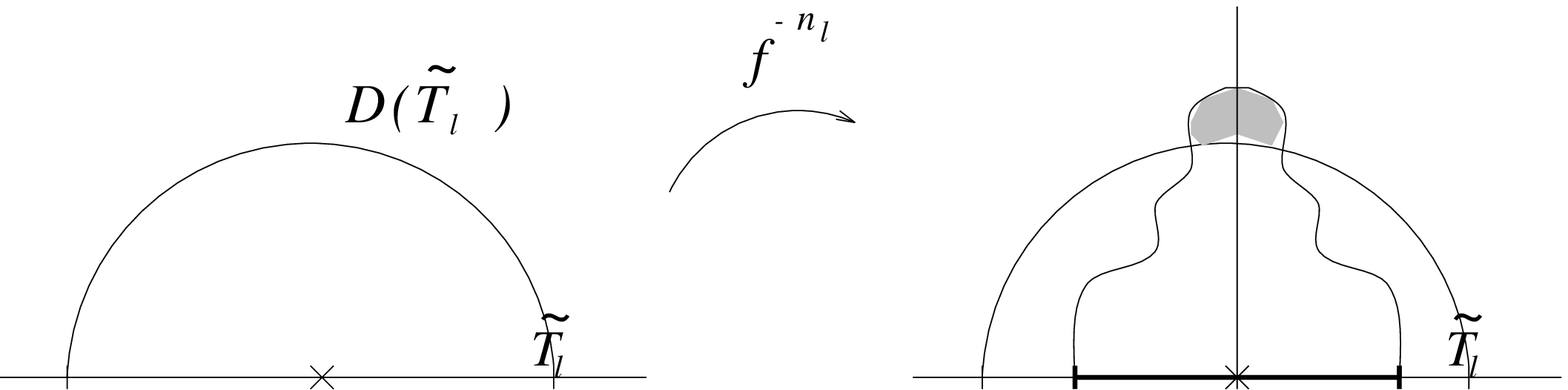}{}{15cm}

Let us consider decomposition (\ref{decompose}). 
The diffeomorphism $\psi_l$ maps some interval $Z^l\supset f  S^l$
onto $T^l$. Hence $\psi_l^{-1}$ has a bounded distortion on $\tl T^l$.
Let $\tl Z^l=\psi_l^{-1} \tl T^l\subset Z^l$.

By bounded geometry, the point $f_l 0$ divides $\tl T^l$
into commensurable parts.  Hence the critical value $f 0=\psi_l^{-1}(f_l0)$
divides $\tl Z^l$ into commensurable parts: Let $A=A(\bar p)$ stand for a bound of the ration of these parts.

By the Schwarz lemma, domain $V=\psi_l^{-1}(D(\tl T^l))$ is contained in $D(\tl Z^l)$.
Hence its pullback $f^{-1}V$ by the quadratic map is contained in a domain
$W=f^{-1}D(\tl Z^l)$ intersecting the real line by $\tl S^l$, and having a bounded distortion
about 0. Hence $\diam W\leq C(\bar p)\tl S^l$, which proves the first statement.

Finally, it follows from \lemref{square root}, that $W\setminus D(\tl T^l)$ is contained in a sector
$Q_{\epsilon}(\tl S^l)$ with an $\epsilon$ depending only on $A=A(\bar p)$ (see Figure 3).

\end{pf}

Let us now give a more precise statement:

\begin{lem}\label{bounds2}
Let $J=J_{-s}$ and $J'=J_{-s'}$
be two returns of the backward orbit (\ref{J-orbit}) to $P^l$,  where $s'=s+tn_l$.
Let $z$ and $z'$ be the corresponding points of the backward orbit of $z_0$. 
Assume $z\in D(\tl T^l)$. 
Then either for some $0\leq i\leq t$, a point $z_{-(s+in_l)}$ $\epsilon$-jumps
 and $|z_{-(s+in_l)}|\leq C |T^l|$,
or  $z_{-s'}\in D_{\theta'}(H')$, where $H'$ is the monotonicity 
interval of $f^{tn_l}$ containing $J'$, and $\theta'=\pi/2-O(\epsilon)$.
\end{lem}

\begin{pf} Assume that the above points do not $\epsilon$-jump. Then by Lemma \ref{bounded1}
they belong to the disk $D(\tl T^l)$. As the map $\psi_l^{-1}$ from \ref{decompose} has a bounded
distortion,  non of the 
points $z_{-m}$ $\delta$-jumps for  $s\leq m\leq s'$, where $\delta=O(\epsilon)$
as $\epsilon\to 0$. Now the claim follows from Lemma \ref{joc}.
\end{pf}

The following lemma will allow us to make an inductive step:

\begin{cor}\label{induction}
Let $J=J_{-(n_l)}$, $J'=J_{-n_{l+1}}$, and $z, z'$ be the corresponding points of the backward
orbit of $z_0$. Assume  $z\in D(\tl T^{l-1})$. Then either there is a good moment $-m\in (-n_l, -n_{l+1})$
 when the point $z_{-m}$ $\epsilon$-jumps and $|z_{-m}|\leq C|T^l|$, or $z'\in D(\tl T^l)$. 
\end{cor}

\begin{pf} Note that by bounded geometry all the moments
 $$-n_l, -(n_l+n_{l-1}), -(n_l+2 n_{l-1}),\ldots, -n_{l+1},$$
when the intervals of (\ref{J-orbit})  return to $P^{l-1}$ before the first return to $P^{l+1}$, are good
(provided the quantifier $\bar K$ is selected sufficiently big).
Hence by \lemref{bounds2} either the first possibility of the claim occurs, or $z'\in D_{\theta'}(L')$, where 
$L'$ is the monotonicity interval of $f^{n_{l+1}-n_l}$ containing $J'$,
and  $\theta'=\pi/2-O(\epsilon)$. As $n_{l+1}-n_l\geq n_l$, 
$L'$ is contained in $S^l$, which is well inside $\tl T^l$. Thus 
$D_{\theta'}(L')\subset D(\tl T^l),$ provided $\epsilon$ is sufficiently small.
\end{pf}

We are ready to carry out the inductive  proof of (\ref{goal}).
Set $j=0$ if $z_0\not\in D(\tl T^0)$.
Otherwise  let $j$ be the smallest level for which
\begin{equation}\label{level}
z_0\in D(\tl T^j).
\end{equation}
By the considerations in the beginning of the proof (if $j=0$) or by \lemref{bounded1} (for $j>0$),
either $z_{-n_j}\in D(\tl T^{j})$, or $z_{-n_j}$ $\epsilon$-jumps. Moreover, in the latter case
$|z_{-n_j}|\leq C |z_0|$, so that (\ref{goal}) follows.

 In the former case we will proceed inductively. Assume that 
either $z_{-n_l}\in D(\tl T^{l-1})$, or $z_{-t}$ $\epsilon$-jumps at some good
moment $-t\geq s$.  If the latter happens, we are done. If the former happens, 
we pass to $l+1$ by Corollary \ref{induction}. \lemref{contraction} is proven.
\subsection{Proof of \lemref{sector} (for bounded combinatorics)}
 By \corref{contr2},  $\diam J(f_k)\leq C|T^k|,$
with a $C=C(\bar p)$. Hence $J(f_k)\subset D(\tl T^l)$, where $l\geq k-N(\bar p)$.
Let $\zeta'\in J(f_k)$, $\zeta=f_k \zeta'$, 
 and $\zeta=\zeta_0,\zeta_{-1}, \ldots, \zeta_{-n}=\zeta'$
be the corresponding backward orbit under iterates of $f_l$. 

By  \lemref{bounded1}, either $\zeta_{-j}$ $\epsilon$-jumps at some moment, or
$\zeta'\in D(\tl T^k)$. If the former happens then $\zeta_{-j}\in D_{\theta}(J_{-j n_l})$,
where $\theta=\theta(\epsilon)>0,$ and $J_{-m}$ are the intervals from \ref{J-orbit}. 
But then by the Schwarz Lemma $\zeta'\in D_{\theta'}(P^k)$ with some $\theta'$ depending
on $\bar p$ only. 
Thus $J(f_k)\subset D_{\theta'}(P^k)\cup D(\tl T^k),$ and we are done.

\begin{rem}
The above proof of the main lemmas for the case of  bounded combinatorics 
 illustrates the ideas involved in treating
the general essentially bounded case.
A complication arises however because of the possibility
that a  jump in the orbit occurs at a "bad" moment when the corresponding
iterate of the periodic interval is not commensurable
with its original size.

\end{rem}

\section{ Essentially Bounded Combinatorics and Geometry}
\label{nest}
Let $f$ be a renormalizable quasi-quadratic map.

We use the standard notations $\beta$ and $\alpha$ for  the fixed 
points of $f$ with positive and negative multipliers correspondingly.
Let $B\equiv B(f)=[\beta,\beta']$, 
$A\equiv A(f)=[\alpha, \alpha']\subset B$.

The map $f$ is called {\it immediately renormalizable} if the interval
$A$ is periodic with period $2$.
If $f$ is not immediately renormalizable, 
let us consider the principal nest
$A\equiv I^{0}\equiv I^0(f)\supset I^{1}\equiv I^1(f)\supset\dots$ of intervals of $f$ (see \cite{L2}).  
It is defined in the following way. Let $t(m)$ be the
first return time of the orbit of 0 back to $I^{m-1}$. Then $I^{m}$
 is defined
as the component of $f^{-t(m)} I^{m-1}$ containing 0. 
Moreover $\cap_m I^{m}=B(Rf)$. 

For $m>1$, let $$g_{m}: \bigcup_{i} I^{m}_i\rightarrow I^{m-1}$$
be the {\it generalized renormalization } of $f$ on the interval $I^{m-1}$,
that is, the first return map restricted onto the intervals intersecting 
the postcritical set (here $I^{m}\equiv I^{m}_0$). Note that 
$g_{m }\equiv f^{t(m)}: I^{m}\rightarrow I^{m-1}$ is unimodal 
with $g_{m} (\partial I^{m})\subset \partial I^{ m-1}$, while
$g_{m}: I^{m}_i\rightarrow I^{m-1}$ is a diffeomorphism for all
$i\not=0$.

Let us consider the following set of levels: 
$$X\equiv X(f)=\{m: t(m)>t(m-1)\}\cup \{0\}=
\{m(0)<m(1)<m(2)<\dots <m(\chi)\}.$$
A level $m=m(k)$ belongs to $X$ iff the return to level $m-1$
is {\it non-central}, that is $g_{m} 0\in I^{m-1}\backslash I^{m}$. 
For such a moment the map $ g_{m+1}$  is essentially different from $g_{m}$ 
(that is not just the restriction of $g_{m}$ to a smaller domain).
 Let us use the notation $h_{k}\equiv g_{m(k)+1}$, $k=1,\ldots \chi$.
The number $\chi=\chi(f)$ 
is called the {\it height } of $f$ ( In the immediately renormalizable
case set $\chi=-1 $).

The nest of intervals 

\begin{equation}\label{cascade}
I^{m(k)+1}\supset I^{m(k)+2}\supset\ldots\supset I^{m(k+1)}
\end{equation}
is called a {\it central cascade}. The {\it length} $l_k$ of the
cascade is defined as $m(k+1)-m(k)$. Note that a cascade of length 1
corresponds to a non-central return to level $m(k)$.

 A cascade \ref{cascade} is called
{\it saddle-node} if $h_k I^{m(k)+1}\not\ni 0$ (see \figref{fig1a}).
 Otherwise it is called
{\it Ulam-Neumann}. For a long saddle-node cascade the map $h_k$
 is combinatorially close
to $z\mapsto z^2+1/4$.  For  a long Ulam-Neumann cascade it is close to 
$z\mapsto z^2-2$. 

\realfig{fig1a}{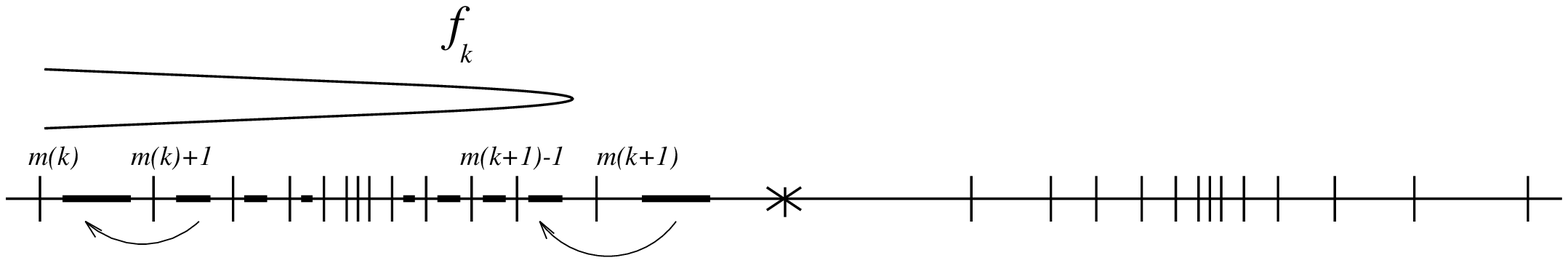}{A long saddle-node cascade}{17cm}

Given a cascade (\ref{cascade}),  let 
\begin{equation}
\label{mark}
K_j^{m(k)+i}\subset
I^{m(k)+i-1}\setminus I^{m(k)+i},\; i=1,\dots,m(k+1)-m(k)-1
\end{equation}
 denote the pull-back of $I_j^{m(k)+1}$ under
 $h_k^{i-1}=g^{i-1}_{m(k)+1}$. Clearly, $K_j^{m(k)+i+1}$
are mapped by $h_k$ onto $K_j^{m(k)+i}$, $i=1,\dots,m(k+1)-m(k)-1$,
while $K_j^{m(k)+1}\equiv I_j^{m(k)+1}$ are mapped onto the whole
$I^{m(k)}$. This family of intervals is called the {\it Markov family} associated 
with the central cascade.

Let $x\in  \omega(0)\cap (I^{m(k)}\setminus I^{m(k)+1})$, 
$h_k x\in I^j\setminus I^{j+1}$. Set 
$$d(x)=\min\{j-m(k), m(k+1)-j\}.$$
This parameter shows how deep the orbit of $x$ lands inside the cascade. 
Let us  now define 
 $d_k$ as the maximum of  $d(x)$ over all
 $x\in \omega(0)\cap (I^{m(k)}\setminus I^{m(k)+1})$. 

Given a saddle-node cascade (\ref{cascade}), let us call all levels
$m(k)+d_k< l < m(k+1)-d_k$ {\it neglectable}. 

Let us now define {\it the essential period} $p_e=p_e(f)$.
 Let $p$ be the period of the periodic interval $J=B(Rf)$.
Let us remove from the orbit $\{f^k J\}_{k=0}^{p-1}$ all intervals 
whose first return to some $I^{m(k)}$ belongs to a neglectable level.  
The essential period is the number of the intervals which are left.

We say that an infinitely renormalizable map $f$
 has {\it essentially bounded combinatorics} if  
$\sup_n p_e(R^n f)< \infty$.  

Let $\sigma(f)=|B(Rf)|/|B(f)|$.
Let us say that $f$ has {\it essentially bounded geometry} if
$\inf_n \sigma(R^n)>0 $. 

\begin{thm}\cite[Theorem D]{L} 
\label{essbound} 
Let $f$ be a quasi-quadratic map of Epstein class. 
There are functions 
 $\delta(p)\geq \epsilon(p)>0$, such that $\delta(p)\to 0$ as
  $p\to\infty$,   with the following properties.
If  $p_e(f)\leq p$  then  $\sigma(f)\geq \epsilon(p)$.
Vice versa, if $p_e(f)\geq p$  then  $\sigma(f)\leq \delta(p)$.
  Thus geometry of $f$ is essentially bounded if and only if its combinatorics is.
\end{thm}

From now on we will work only with maps having essentially bounded combinatorics,
and $\bar p$ will stand for a bound of the essential period. 
By the {\it gaps} $G^m_j$ of level $m$ we mean  the components of $I^{m-1}\setminus \cup I^m_j$.
We say that a level $m$ is {\it deep inside the cascade} if $m(k)+\bar p\leq m\leq m(k+1)-\bar p$. 
Let us finish this section with a lemma on geometry of maps with essentially bounded combinatorics.

\begin{lem}\cite [Lemma 17]{L}
\label{commens}
Let  $f$ be a quasi-quadratic map  with essentially bounded combinatorics. Then for any $m$,
 the non-central intervals $I^m_i,\; i\not=0,$ and  the gaps $G^m_j$ of level $m$ 
are $C(\bar p)$-commensurable with $I^{m-1}\setminus I^m$. Moreover, this is also true for
the central interval $I^m_0$, provided $m$ is not deep inside the cascade.
\end{lem}

Note that the last statement  of the lemma is definitely false when $m$ is deep inside the
cascade: then $I^m$ occupies almost the whole of $I^{m-1}$. So we observe commensurable
intervals in the beginning and in the end of the cascade, but not in the middle.
  This is the saddle-node  phenomenon which is in the focus of this work.

\section{Saddle-Node Cascades}
\label{cascades}
 
Let $f\in \EE_\lambda$ be a map of Epstein class.
\comm {In this section we fix the level $i$ and work with the map
$F\equiv f_i$. Let $I^0\supset I^1\supset\dots$ be the 
principal nest of $F$, and let $g_m$ be the generalized 
renormalization of $F$ on the interval $I^{m-1}$.

 Let the intervals
\begin{equation}
\label{cscd}
I^{m(k)+1}\supset I^{m(k)}\supset \dots\supset I^{m(k+1)} 
\end{equation}
 of the principal nest form  a  central cascade. 
 
Following the notations of  \secref{nest} we set $h_{k}\equiv g_{m(k)+1}$.
We assume the central cascade (\ref{cscd}) to be saddle-node, in which 
case one has $h_k I^{m(k)+1}\not\ni 0$.
the {\it length of the cascade} is $l=m(k+1)-(m(k)+1)$. }

Let us note first for a long saddle-node cascade \ref{cascade}, the  map
 $h_{k}:I^{m(k)+1}\rightarrow I^{m(k)}$ is  a 
small perturbation of a map with a parabolic fixed point.

\begin{lem}\cite{L} \label{pert}
Let $h_k$ be a sequence of maps of Epstein class $\EE_{\lambda}$
having saddle-node cascades of length $l_k\to\infty$. Then
any limit point $f: I'\rightarrow I$ 
of this sequense (in the Caratheodory topology) 
has on the real line topological  type of $z\mapsto z^2+1/4$, and thus
has a parabolic fixed point. 
\end{lem}

\begin{pf} 

It takes $l_k$ iterates for the critical point to escape  $I^{m(k)+1}$
under iterates of $h_k$. Hence the critical point does not 
escape $I'$ under iterates of $f$. By the kneeding theory \cite{MT}
$f$ has on the real line topological type of $z^2+c$ with $-2\leq c\leq 1/4$.
Since small perturbations of $f$ have escaping critical point,
the choice for $c$  boils down to only two boundary 
parameter values, $1/4$ and $-2$.
Since the cascades of $h_k$ are of saddle-node type, $f I'\not\ni 0$, which
rules out  $c=-2$.  

\end{pf}

\begin{rem}
Thus  the plane dynamics of $h_k$ with a
 long saddle node cascade 
resembles the dynamics of a map with a parabolic fixed point: the
orbits follow horocycles (cf. \figref{fig6}).
\realfig{fig6}{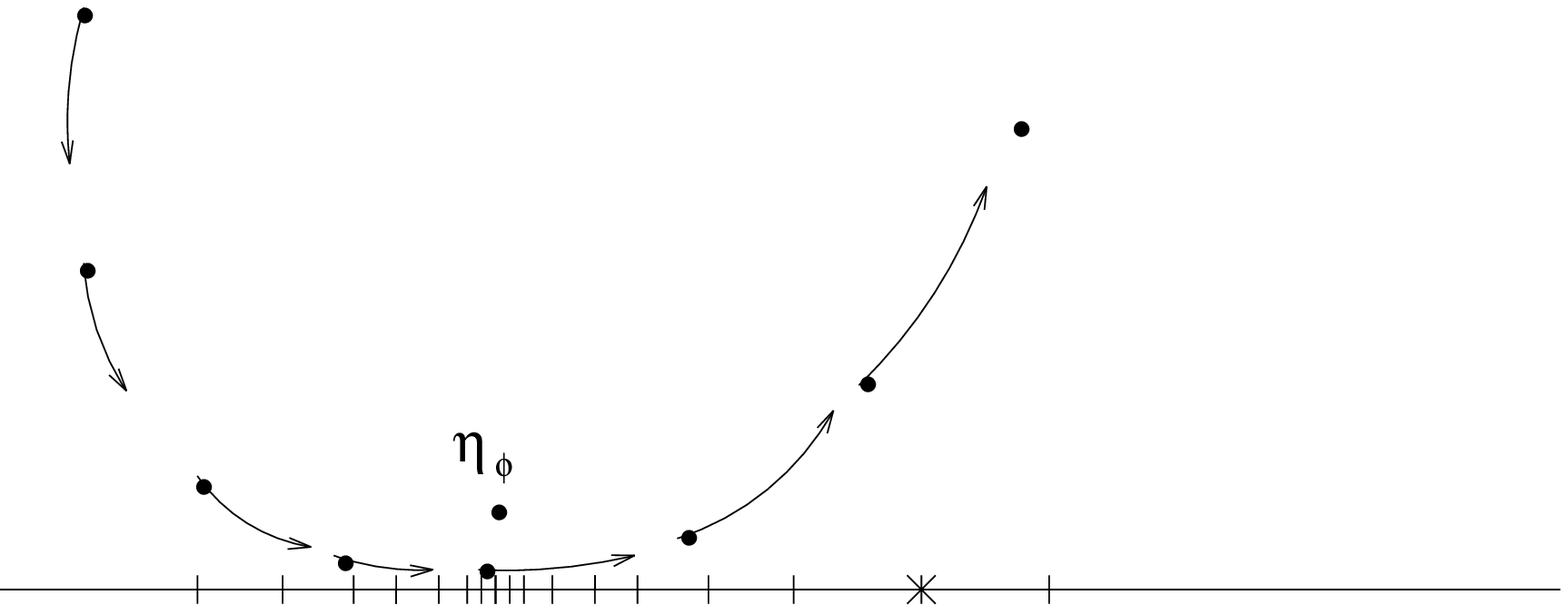}{The backward trajectory of a point corresponding
to a saddle-node cascade}{13 cm}
\end{rem}
 
\comm{For a level $k>i$ we consider the backward orbit
\begin{equation}
\label{J-orbit-1}
J_0\equiv P_0^k, J_{-1}\equiv P^k_{n_k-1},\dots,J_{-n_k-1}\equiv P_1^k.
\end{equation}}

\begin{lem}
\label{casc1}  Let us consider a saddle-node cascade \ref{cascade} generated
by a return map $h_k$. 
Let us aslo consider a backward orbit of an interval $E\subset I^{m(k)}\setminus I^{m(k)+1}$ under 
iterates of $h_k$:
 $$E\equiv E_0,\;
  E_{-1}\subset I^{m(k)+1}\setminus I^{m(k)+2} ,\dots, 
E_{-j}\equiv E'\subset I^{m(k)+j}\setminus I^{m(k)+j+1},$$ where
$ m(k)+j+i\leq m(k+1)$.  
Let $z=z_0,z_{-1},z_{-2},\dots,z_{-j}=z'$ be the corresponding  
 backward orbit of a point  $z\in D(I^{m(k)})$. If the length of the
cascade is sufficiently big, then either 
 $z'\in D(I^{m(k)})$, or  $\ang{z'}{J'}>\epsilon$ and
$\dist (z',J')\leq C(\bar{p})|I^{m(k)}|$.
\end{lem}

\begin{pf} To be definite, let us assume that the intervals $E_{-i}$ lie on the
left of 0 (see Figure 4). Without loss of generality, we can assume that
$z\in \Bbb H$.
Let $\phi=h_k^{-1}$ be the inverse branch of $h_k$ for which
$\phi E_{-i}=E_{-(i+1)}$. As  $\phi$ is orientation preserving on 
$(-\infty, h_k 0]$, it maps the  upper half-plane $\Bbb H$
into itself: 
  $\phi(\Bbb H)\subset\{z=r^{ei\theta}|\;r>0,\;\pi>\theta >\pi/2\}$.

\comm{ univalent maps $\phi$ with the above property form a normal family 
in $\Bbb H$ by Montel's theorem.}

By \lemref{pert}, if the cascade \ref{cascade} is sufficiently long,
the map $\phi$ has an attracting fixed point
 $\eta_{\phi}\in \Bbb H\cap D(I^{m(k)+2}) $ (which is a perturbation
of the parabolic point for some map of type $z^2+1/4$). 
By the Denjoy-Wolf Theorem, 
$\phi^{n}(\zeta)\underset{n\rightarrow\infty}
\longrightarrow \eta_{\phi}$ for any $\zeta\in\Bbb H$, uniformly on compact
subsets of $\Bbb H$.
Thus 
for a given compact set $K\Subset \Bbb H$, there exists $N=N(K,\phi)$ such that 
$ \phi^N(K)\subset D(I^{m(k)+1})$. By a normality argument,
the choice of $N$ is actually independent of a particular $\phi$ under 
consideration.

Suppose $z_{-r}\notin D(I^{m(k)})$.
By \lemref{square root} the set 
$K= D(I^{m(k)})\setminus \phi(D(I^{m(k)}))\cap \Bbb H$ is compactly contained 
in $\Bbb H$, and $\diam K\leq C|I^{m(k)}|$. For $N$ as above we have 
$z'\in \cup_{i=0}^{N-1}\phi^i(K)\cup D(I^{m(k)})$ and the lemma is 
proved.

\end{pf}

\comm{
\begin{lem}
\label{casc2}
Given a central cascade (\ref{cscd}), let  $L_{s}$ be a component of $I^{m(k+1)-s}\setminus I^{m(k+1)-(s-1)}$, $L_{s}\equiv
 g_{m(k)+1}^s(I^{m(k+1)})$, $0\leq s\leq l$. 
Then the pull-back $\phi^s : D_\theta(L_s)\to D_\theta(I^{m(k+1})$
has distortion bounded by a constant depending on $\theta$ and $\bar p$
only.
\end{lem}
\realfig{fig4}{fig4.eps}{}{12 cm}
\begin{pf}
The geodesic neighbourhood $D_\theta(L_{s})$ is contained inside $D_\theta(J_{s+1}\cup
J_{s}\cup J_{s-1})$. The annulus $D_\theta(J_{s+1}\cup
J_{s}\cup J_{s-1})\setminus D_\theta(J_{s})$ has a definite modulus 
depending on $\theta$ (cf. \figref{fig4}).
By  Koebe Theorem the mapping $\phi^{s-1}:D_\theta(L_{s})\rightarrow
 D_\theta(L_1)$ has bounded distortion, and the lemma follows.
\end{pf}
From the above lemma we derive a corollary:
\begin{cor}
\label{casc3}
Let $J=J^1_0=J_{-i},J^1_{-1},\dots, J^{1}_{-j}=J'$ be the 
returns of the orbit (\ref{J-orbit-1}) to $I^{m(k)}$, between two
consecutive returns $J$ and $J'$ to $I^{m(k+1)}$.
Let $z=z_0,z_{-1},z_{-2},\dots,z_{-j}=z'$ be the corresponding
points of the backward orbit of $z$.
Let $z\in D_\theta(I^{m(k)})$. Then $z'\in D_{\theta'}(I^{m(k+1)})$ with
 $\theta'=\theta'(\theta)$. Moreover, $\ang{z'}{J'}=O(\ang{z}{J})$.
\end{cor}
\begin{pf}
To control the distortion of the pull-back, we use the Markov scheme
(\ref{mark}), associated with the cascade.

Let us decompose the orbit  $J^1_{-1},\dots, J^{1}_{-j}=J'$
into the pieces $J^1_{-t(r)}, J^1_{-t(r)-1},\dots,   
J^1_{-t(r)-s}$,
such that $J^1_{-t(r)}\subset K_j^{m(k)+2},\;
J^1_{-t(r)-1}\subset K_j^{m(k)+3},\dots,
J^1_{-t(r)-i}\subset K_j^{m(k)+i+2},$ for $i<s$,
and $J^1_{-t(r)-s}\subset I^{m(k)}\setminus I^{m(k+1)}$.

There are finitely many  such pieces by essentially bounded 
combinatorics. The distortion control for the corresponiding
piece of the orbit
of $z$ is provided by \lemref{casc2}.
\end{pf}

} 

\section{Proofs of the Main Lemmas}
\label{inductive}

\subsection{ Proof of \lemref{contraction}}
Let us start with a little lemma:
\begin{lem}
\label{divides}
Let $f\in \EE_\lambda$ be a map  of Epstein class without attracting
 fixed points.
Then both components of $B\setminus A$   contain an $f$-preimage of $0$
which divides them into $C(\lambda)$-commensurable parts.
\end{lem}
\begin{pf}
The interval $[\alpha,\beta']$ is mapped by $f$ onto $[\beta, \alpha]\ni 0$.
Denote by $\eta=f^{-1}(0)\cap [\alpha,\beta']$. Under our assumption
this point is clearly different from $\alpha$ and $\beta'$.
As the space of maps of Epstein class $\EE_\lambda$ with no attracting
fixed points is compact,
$\eta$ divides $[\alpha,\beta']$ into $C(\lambda)$-commensurable parts.
The analogous statement is certainly true for the symmetric 
point $\eta'\in [\beta,\alpha']$.
\end{pf}

\label{ess}
 As in \secref{bounded}, let us fix a level $\tau$, let $n=n_\tau$,
 and set  
\begin{equation}
\label{J-orbit-2}
J_0\equiv P^\tau, J_1\equiv P^\tau_{n-1},\dots, J_{-(n-1)}\equiv P^\tau_1.
\end{equation} 

For any point $z\in \C{T^\tau}$ with $\dist(z,J_0)>|J_0|$, we denote by
\begin{equation} \label{z-orbit}
z\equiv z_0,z_{-1},z_{-2},\dots,z_{-(n-1)}
\end{equation}
the backward orbit of $z$ corresponding to the orbit ($\ref{J-orbit-2}$).
We should  prove that
\begin{equation}
\label{goal-2}
{\dist(z_{-(n-1)}, J_{-(n-1)})\over |J_{-(n-1)}|}\leq C(\bar p) 
  {\dist (z_0, J_0)\over |J_0|}.
\end{equation}

Let $A=A(f)$ and $B=B(f)$ be the intervals defined in \secref{nest},
$H_s(x)$ be the monotonicity intervals as defined in \secref{jumping}. 

\begin{lem}[First return to $A$]\label{beginning}
Let $ J_{-s}$ be the first return of the orbit (\ref{J-orbit-2}) to the 
interval $A$.
 There is an $\eps=\eps(\bar p)>0$ such that  either 
$z_{-s}\in D(B),$ 
or there is a moment $-i\in [ -s, 0]$ when the backward orbit (\ref{z-orbit})
$\epsilon$-jumps: $\ang{z_{-i}}{J_{-i}}>\epsilon$ and moreover
\begin{equation}\label{sqrt-2} 
{\dist(z_{-i}, J_{-i})\over |J_{-i}|}\leq C_0(\bar p) 
  {\dist (z_0, J_0)\over |J_0|}.
\end{equation} 
\end{lem}

\begin{pf}  By  definition of the essential period,
 $s\leq p_e(f)\leq \bar p$. Hence (\ref{sqrt-2}) holds
for  all $ i=0,1,\ldots,  s$ by the first
statement of \lemref{square root2}.

Further, by \lemref{divides}
each component of $B\setminus A$ contains a preimage of $0$, and it divides 
$B$ into $ K(\bar p)$-commensurable intervals. Hence, the monotonicity
interval $H=H_s(J_{-s})$ is well inside of $B$,
and the conclusion follows from  \lemref{joc}. 
\end{pf}

If the second possibility of \lemref{beginning} occurs then
  (\ref{goal-2}) follows from  
\lemref{goodangle}. If the first one happens, we proceed inductively along the 
principal nest. Namely, in the following  series of lemmas 
we will show that 
the backward $z$-orbit (\ref{z-orbit}) either $\eps$-jumps at some
good moment, or follows the backward $J$-orbit (\ref{J-orbit-2}) with at most
one level delay.

In the  following lemmas we work with a fixed renormalization level $l$
 and skip  index $l$ in the
notations: $f\equiv f_l\equiv  R^l (f_0)$,
  $A\equiv A^l, B\equiv B^l$. We will use notations of \secref{nest}
for different combinatorial objects. 

\comm{
Recall that $g_m=\cup I^m_i\rightarrow I^{m-1}$ is the generalized 
renormalization on the interval $I^{m-1}$. We mark the levels with
noncentral returns $A\supset I^{m(1)}\supset I^{m(2)}\supset\dots\supset
 I^{m(\chi)}\supset B'$.

 The levels between two noncentral returns form a central
cascade
\begin{equation}
I^{m(k)}\supset I^{m(k)+1}\supset\dots\supset I^{m(k+1)} 
\end{equation}
 generated by the return map $f_k\equiv g_{m(k)+1},\;
 f_{k}(0)\in I^{m(k)}\setminus I^{m(k)+1}$.                end comm }

\begin{lem}[Further returns to $A$]
\label{ind1}
Let $E=E_0, E_{-1}, \dots , E_{-s}=E'$ be the consecutive returns of  
the backward orbit (\ref{J-orbit-2}) to $B$, between
two consecutive returns to $A$. Let 
 $\zeta=\zeta_0, \zeta_{-1},\dots \zeta_{-s}=\zeta'$
 be the corresponding  points of the backward 
orbit (\ref{z-orbit}). Assume $\zeta\in D(B)$.
Then either $\zeta'\in D(B)$, or
 $\ang{\zeta_{-i}}{E_{-i}}>\epsilon(\bar p)>0$ and
 $\dist(\zeta_{-i}, E_{-i})\leq C(\bar{p})|B|$ 
for some $0\leq i\leq s$.
\end{lem}

\begin{pf} Take an $\epsilon>0$.
 By definition of the essential period,  
$s\leq \bar p$. 
By the Schwarz lemma and \lemref{square root},  $\zeta_{-i}\in D_\sigma(B)$ for  $i=0,1\dots s$, with $\sigma=\sigma(\bar p)$.
 Hence $ \dist(\zeta_{-i}, E_{-i})\leq C(\bar p)|B|.$
 If for some $i\in [0,s]$, 
$ \ang{\zeta_{-i}}{E_{-i}}>\eps,$
we are done. 

By \lemref{divides} each component of $B\setminus A$ contains
  an $f$-preimage of 0 which divides $B$ into $K$-commensurable
 intervals,
with $ K =K(\bar p)$. Hence the monotonicity interval of $f$,
 $H= H_s(E_{-s} )$, 
is well inside of $B$. As $f: B\rightarrow B$  has an extension of Epstein
 class $\EE_{\lambda}$ (see \S 2.6), we can apply 
 \lemref{joc-2}. It follows that if none of the points $\zeta_{-i}$ $\epsilon$-jumps, then $\zeta_{-i}\in D_{\theta}(H)$, $0\geq -i\geq -s$, with $\theta=\pi/2-O(\eps)$. Thus $\zeta_{-s}\in D(B)$ for sufficiently small $\epsilon<\eps(\bar p)$,
 and the proof is completed. 
\end{pf}

We say that a point/interval is deep inside of the cascade (\ref{cascade})
if it belongs
to $I^{m(k)+\bar p}\backslash I^{m(k+1)-\bar p}$. (In the case of essentially
bounded combinatorics  such a cascade must be of
 saddle node type). Recall that a moment $-i$ is called ``good" if the 
interval $J_{-i}$ is commensurable with $J_0$. Because of  the essentially
bounded geometry, this happens, e.g., when for some $k$, the interval
$J_{-i}$ lies in $I^{m(k)}\setminus I^{m(k+1)}$ but is not deep inside the 
corresponding cascade.

  \begin{lem}[First return to $I^{m(1)}$]
\label{ind2}
 Assume that $f$ is not immediately renormalizable.
Let $E\equiv  E_0, E_{-1}, \dots,$ $ E_{-s}\equiv E'$ 
be the consecutive returns of the 
backward orbit
 (\ref{J-orbit-2}) to $A$ until the first return to $I^{m(1)}$.
Let $ \zeta\in \C{A} \cap D(B)$, and  let 
$\zeta\equiv \zeta_0, \zeta_{-1}\dots \zeta_{-s}\equiv \zeta'$ be the corresponding points in the backward orbit of $\zeta_0$.  
Then either $\zeta'\in D(A)$, 
or $\ang{\zeta_{-i}}{E_{-i}}>\epsilon(\bar p)>0$ and
 $\dist(\zeta_{-i}, E_{-i})\leq C(\bar p)|B|$ 
at some good moment $0\geq -i\geq -s$.
\end{lem}

\begin{pf} Let $H=H_s(E_{-s})$.

  As $f$ is not immediately renormalizable,
we have the interval
$I^{1}=[p,p']$. Let   $p$ be chosen on the same side of $0$
as $\alpha$. 
 Then $f^2[\alpha, p]\supset[\alpha, \alpha']$.
  Denote by $\eta$
the $f^2$-preimage of $0$ in $[\alpha, p]$. Since $f$ is quadratic up to
bounded distortion, the map 
$f^2|_{[\alpha, p]}$ is quasi-symmetric (that is,
maps commensurable adjacent intervals onto commensurable ones). It follows
that   $\eta$
divides $[\alpha, p]$, and hence $A$,
 into $ K=K(\bar p)$-commensurable parts. 
Hence $H\subset [\eta, \eta']$ is well inside $A$. 

By  \lemref{joc-2} and  \lemref{ind1}, either $\zeta'\in D_{\theta}(H)$ with
 $\theta=\pi/2-O(\eps)$, or there is a moment $i\leq s$
such that 
\begin{equation}
\label{ind2-1}
\ang{\zeta_{-i}}{E_{-i}}> \epsilon \quad {\rm and} \quad
\dist(\zeta_{-i}, E_{-i}) \leq C(\bar p)|B|.    
\end{equation}

In the former case we are done as $D_{\theta}(H)\subset D(A)$ for sufficiently 
small $\eps$. 

Let the latter case occur. Then we are done if the moment  $-i$ is good.
 Otherwise  $E_{-i}$ is deep inside the cascade
$A=I^0\supset I^1\supset \dots \supset I^{m(1)}$.
Consider the largest $r$ such that
$E_{-(i+q)}\subset I^{t+q-1}\setminus I^{t+q}$
 for all $0\leq q\leq r$. Note that by
essentially bounded combinatorics, the moment $-j=-(i+r)$ has to be good. 
By \lemref{casc1},
either (\ref{ind2-1}) occurs for $\zeta_{-j}$, and we are done, or
 $\zeta_{-j}\in D(A)$.

In the latter case let $\tl K\subset I^{m(1)-1}\setminus I^{m(1)}$ be
the interval containing $E_{-(s-1)}$ which is homeomorphically mapped
under $h_1^{s-1-j}$  onto $A$ 
(to see that such an interval exists, 
consider the  Markov scheme described in \S 5). By the Schwarz lemma
 $\zeta_{-(s-1)}\in D(\tl K) \subset D(A)$.
 Now the claim follows from  \lemref{square root}.
 \end{pf}

Now we are in a position to proceed inductively along the principal nest:
Note that the assumption of the following lemma is checked for $k=1$
in \lemref{ind2}.  

\begin{lem}[Further returns to $I^{m(k)}$]
\label{ind3}
Let $E$ and $E'$ be two consecutive returns of the
backward orbit
(\ref{J-orbit-2}) to the interval $I^{m(k)}$.
Let $\zeta$ and $\zeta'$ be the corresponding points of the backward orbit
of $z_0$.
Assume that $\zeta\in D(I^{m(k-1)})$. Then, either $\zeta'\in D(I^{m(k)})$, or
$\ang{\zeta'}{E'}>\epsilon(\bar p)>0$, and 
$\dist(\zeta',E')<C(\bar p)|I^{m(k-1)}|.$
\end{lem}
\begin{pf}
Denote by $\tl E$ the last interval in the backward orbit (\ref{J-orbit-2})
 between $E$ and $E'$, which
visits $I^{m(k-1)}$ before returning to $I^{m(k)}$. Then $h_kE'=\tl E$ and
$h_k^{\circ j} \tl E=E$ for an appropriate $j$.

The Markov scheme (\ref{mark}) provides us with an interval 
$\tl K\subset I^{m(k-1)}\setminus I^{m(k)}$ containing $\tl E$ which is homeomorphically mapped 
under $h_k^{\circ j}$ onto $I^{m(k-1)}$. By essentially bounded geometry
and distortion control along the cascade\marginpar{where is it?!},
$\tl K$ is well inside  $I^{m(k)-1}\setminus I^{m(k)}$,
 and the critical value of $h_k$
divides $\tl K$ into commensurable parts.

Let $K'\supset E'$ be the pull-back of $\tl K$ by $h_k|I^{m(k)}$.
It follows that  $K'$ is contained well inside $ I^{m(k)}$. 

Let $\tl \zeta=h_k \zeta'$ be the point  of the orbit (\ref{z-orbit}) corresponding to $\tl E$.
By the Schwarz lemma, $\tl \zeta\in D(\tl K)$. 
 By the previous remarks and  \lemref{square root},
  $\zeta'\in D(I^{m(k)})$, or
$\ang{\zeta'}{E'}>\epsilon(\bar p)$
 and $\dist(\zeta',E')<C(\bar p)|I^{m(k-1)}|.$
\end{pf}

\lemref{ind3} is not enough for making inductive step since the jump can
occur at a bad moment. The following lemma takes care of this possibility
in the way similar to \lemref{ind2}.

\begin{lem}[First return to $I^{m(k+1)},\; k\geq 1$]
\label{ind4}
 Let  
 $E\equiv E_0, E_{-1}, \dots, E_{-s}\equiv E' $ be the consecutive 
returns of the orbit (\ref{J-orbit-2}) to $I^{m(k)}$ until the
first return to $I^{m(k+1)}$. Let
$\zeta\equiv \zeta_0, \zeta_{-1}\dots,\zeta_{-s}\equiv \zeta' $ be the
 corresponding points
in the backward orbit of $\zeta$. Assume that
 $\zeta_{-1}\in \C{I^{m(k+1)}} \cap D(I^{m(k)}) $.
  Then either $\zeta'\in D(I^{m(k)})$, or $\ang{\zeta_{-i}}{E_{-i}}>
\epsilon(\bar p)>0 $
and $\dist(\zeta_{-i},E_{-i})<C(\bar p)|I^{m(k)}|$
 at some good moment $-1\geq -i\geq -s$.
\end{lem}
\begin{pf}
Let $H\supset E'$ be the maximal interval on which $f_k^{\circ s}$ is monotone.
Note, that both components of $I^{m(k)}\setminus I^{m(k)+1}$ contain
 pre-critical 
values of $h_k$, which divide $I^{m(k)}$ into $K(\bar p)$-commensurable
parts. Hence, $H$ is well inside of $I^{m(k)}$.

By  \lemref{joc-2}, either $\zeta'\in D_{\theta}(H)$ with
 $\theta=\pi/2-O(\eps)$, or there is a moment $1\leq i\leq s$
such that 
\begin{equation}
\label{ind2-2}
\ang{\zeta_{-i}}{E_{-i}}> \epsilon \quad {\rm and} \quad
\dist(\zeta_{-i}, E_{-i}) \leq C(\bar p)|I^{m(k)}|.    
\end{equation}

In the former case we are done as $D_{\theta}(H)\subset D(I^{m(k)})$ if 
$\eps$ is sufficiently small. 

Let the latter case occur. Then we are done if the moment  $-i$ is good.
 Otherwise  $E_{-i}$ is deep inside the cascade
$I^{m(k)}\supset I^{m(k)+1}\supset \dots \supset I^{m(k+1)}$.
 Consider the largest $r$ such that
$E_{-(i+q)}\subset I^{m(k)+t+q-1}\setminus I^{m(k)+t+q}$ for all $q\leq r$.
 Note that by
essentially bounded combinatorics, the moment $-j=-(i+r)$ has to be good. 
By \lemref{casc1},
either (\ref{ind2-2}) occurs for $\zeta_{-j}$, and we are done, or
 $\zeta_{-j}\in D(I^{m(k)})$.

In the latter case, the Markov scheme (\ref{mark}) provides
us with an interval $\tl K\subset I^{m(k+1)-1}\setminus I^{m(k+1)}$
containing $E_{-(s-1)}$ which is mapped homeomorphically onto $I^{m(k)}$ by
$h_k^{s-1-j}$. By the Schwarz Lemma $\zeta_{-(s-1)}\in D(\tl K)\subset
D(I^{m(k)})$. The claim now follows from \lemref{square root}.
\end{pf}
The following lemma will allow us to pass to the nest renormalization level.
It is similar to \lemref{ind1} except that  we deal with a map
of Epstein class rather than a quadratic map. Let us restore now 
label $l$ for the renormalization level.

\begin{lem}[To the next  renormalization level: period $>2$ case] 
\label{ind5}
Assume$ \;$ that $f_l$ is not immediately renormalizable.
\comm{Let $E= E_0,E_{-1},\dots,$  $E_{-s}=E'$ be the returns
of the backward orbit (\ref{J-orbit-2}) to $B^{l+1}$ until the first return
 to $A^{l+1}$. 
Let $\zeta=\zeta_0,\dots, \zeta_{-s}=\zeta'$ be the
 corresponding points of the backward 
orbit (\ref{z-orbit}). 
Assume  that  $\zeta\in D(I^{m(\chi)})$, where
$I^{m(\chi)-1}$ is the last non-central level of the principal nest of $f_l$.
 Then  either $\zeta'\in D(B^{l+1})$, or
 $\ang{\zeta_{-i}}{E_{-i}}>\epsilon(\bar p)>0$ and 
 $\dist(\zeta_{-i}, E_{-i})\leq C(\bar{p})|B^{l+1}|$
for some $0\leq i\leq s$. end comm}
Let $E=E_{-1},\dots,E_{-r}=E',\dots,E_{-(r+s)}=E''$ be the returns
of the backward orbit (\ref{J-orbit-2}) to $B^{l+1}$, and let $E'$, $E''$
  be two consecutive
returns to $A^{l+1}$.
Let $\zeta=\zeta_{-1},\dots,\zeta',\dots, \zeta_{-(r+s)}=\zeta''$ be the
 corresponding
 points of the backward 
orbit (\ref{z-orbit}), and suppose $\zeta\in D(I^{m(\chi-1)})$, where
$\chi=\chi(f_l)$ is the height of $f_l$.
 Then  either $\zeta''\in D(B^{l+1})$, or
 $\ang{\zeta_{-i}}{E_{-i}}>\epsilon(\bar p)>0$ and 
 $\dist(\zeta_{-i}, E_{-i})\leq C(\bar{p})|B^{l+1}|$
for some $1\leq i\leq r+s$.
 Moreover, all these moments are good.
\end{lem}

\begin{pf}
 First, $r+s\leq 2\bar p$ by definition of the essential period
 $\bar p$, and the last statement follows.

 By \lemref{ind3}, either $\ang{\zeta_{-2}}{E_{-2}}>\eps$, $\dist(\zeta_{-2},
E_{-2})\leq C(\bar p)|B^{l+1}|$, or $\zeta_{-2}\in
D(I^{m(\chi)})$.
 
By the Schwarz lemma and \lemref{square root}, 
if $\zeta_{-i}\in D(I^{m(\chi)})$, then  either 
$\zeta_{-(i+1)}\in D(I^{m(\chi)})$, or
$\dist(\zeta_{-(i+1)}, E_{-(i+1)})\leq C(\bar p)|B^{l+1}|$ and 
 $\ang{\zeta_{-(i+1)}}{E_{-(i+1)}}>\epsilon(\bar p)>0$.
 In the latter case we are done.

If the former case occurs for all $i<r+s$ then
by \lemref{joc-2},
 $\zeta''\in D_{\theta}(H) $, where  $H= H_{r+s-1}(E'', f_{l+1})$ and 
 $\theta=\pi/2-O(\eps)$.       By \lemref{divides}, $H$ is well inside
 $B^{l+1}$, and hence
$D_{\theta}(H)\subset D(B^{l+1})$ 
for sufficiently small $\epsilon>0$. 
\end{pf}

Our last lemma takes care of the case when the map $f_l$ is immediately 
renormalizable.

\begin{lem}[To the next renormalization level: period 2 case]
\label{immed}
Assume that $f_l$ is immediately renormalizable, so that $A^l=B^{l+1}$. 
Let $E\subset B^{l+1}$, 
$E\equiv E_0,E_{-1},\dots,$$E_{-s}$$\equiv E'$ be the consecutive returns
 of the backward orbit (\ref{J-orbit-2}) to $B^{l}$, 
until the first return to
$A^{l+1}$.
 \comm{\marginpar{Reason for changes: Lemma 4.2 is not applicable to $f_{l+1}$
                       as it is not extended over $B^l$}}
Let
 $\zeta\equiv \zeta_0,\dots,\zeta_{-s}\equiv\zeta'$
be the corresponding points of the backward orbit (\ref{z-orbit}).
Assume also that  $\zeta\in \C{A^{l}}\cap D(B^l)$.
Then either $\zeta'\in D(B^{l+1})$, or 
$$\ang{\zeta_{-i}}{E_{-i}}>\eps\quad and \quad 
 \dist(\zeta_{-i},E_{-i})<C(\bar p)|B'|$$ for some $0\geq -i\geq -s$.
Moreover, all these moments are good.
\end{lem}
\begin{pf}
By essentially bounded combinatorics, $s\leq 2\bar p$ which yields the 
last statement.

Further, by \lemref{divides},
the monotonicity interval $H_s(E_{-s}, f_l)$ is contained well inside of
 $B^{l+1}$, and the claim follows from  \lemref{joc-2}.
 \end{pf}

\comm{\marginpar{would be good to notice in the statements of lemmas
that we actually finish "inside" the
monotonicity intervals (when we don't jump)}}

Let us now summarize the above information. 
When $f_{\tau-1}$ is immediately renormalizable, set $V_\tau=B^{\tau -1}$.
Otherwise let $ V_\tau= I^{m(\chi-1)-1}(f_{\tau-1})$ where $\chi=\chi(f_{\tau-1})$ is
the height of $f_{\tau-1}$.



\begin{lem}\label{summary} Let $f_\tau=R^\tau f$. 
Let  us consider the backward orbit (\ref{J-orbit-2}) of an interval $J$
and the corresponding orbit (\ref{z-orbit}) of a point $z$. Then
there exist $\eps=\eps(\bar p)>0$ 
  such that
either one of the points $z_{-s}$
$\eps$-jumps at some good moment, or  $z_{-(n-1)}\in D(V_\tau)$.
\end{lem}

\noindent{\bf Proof of Lemma 3.1.} If the 
 former possibility of \lemref{summary} occurs than \lemref{goodangle}
yields (\ref{goal-2}). In the latter possibility happens then
$${    \dist(z_{-(n-1)}, J_{-(n-1)})
\over |J_{-(n-1)}|   }    \leq C(\bar p)$$
by essentially bounded geometry, and we are done again.

Lemma 3.1 is proved.  \QED

\comm{
Let us now put everything together.
Set $l=0$ if $z_0\notin D(B^0)$; otherwise let $l$ be the smallest level
for which $z_0\in D(B^l)$. 
Let us mark by $J'$ the first return of  the backward orbit
(\ref{J-orbit-2}) to $A^l$,
 and let $z'$ be the corresponding point 
of the orbit (\ref{z-orbit}). By \lemref{beginning} (if $l=0$) or
by \lemref{ind1}, either
 one of the points $z_0, z_{-1},\dots,z_{-t}\equiv z'$
$\epsilon$-jumps, or $z'\in D(B^l)$. In the former case
(\ref{goal-2}) follows from \lemref{goodangle}. Otherwise
 we will go further along the principal nest.

Assume first that $f_l$ is not immediately renormalizable.
Denote by $E_0,E_{-1},\dots,$ $E_{-s}\equiv E',\dots,$ $E_{-2s}\equiv E''$
 the consecutive
returns of the orbit (\ref{J-orbit-2}) to the element $I^{m(k)}$ of
the principal nest of $F_l$, where $E',$ $E''$ mark the two
consecutive returns to $I^{m(k+1)}$. Let
 $\zeta_0,\dots,\zeta_{-s}\equiv\zeta',\dots,\zeta''$ be the
 corresponding points in the orbit (\ref{z-orbit}).

Suppose, $\zeta_0\in D(I^{m(k)})$ if $k\geq 1$, or $\zeta_0\in D(B_l)$
otherwise.
By Lemmas \ref{ind2}, \ref{ind3}, either the orbit (\ref{z-orbit})
$\epsilon$-jumps at a good moment $\zeta_{-i}$ when $E_{-i}$ is commensurable
with $J_0$ for some $i\leq s$, or $\zeta'\in D(I^{m(k)})$.
If the former happens, the claim of (\ref{goal-2}) follows
from \lemref{goodangle}.
Let the latter happen. Then by \lemref{ind4}, $\zeta''\in D(I^{m(k+1)}).$

By induction on $k$, either the orbit (\ref{z-orbit}) $\eps$-jumps
at one of the good moments described above, or $\hat \zeta\in D(I^{m(\chi)})$,
where $I^{m(\chi)}$ is the last non-central element of the principal nest
of $F_l$, and $\hat \zeta$ is the point in the orbit (\ref{z-orbit}), which 
corresponds to the second return of (\ref{J-orbit-2}) to $I^{m(\chi)}$.

Now we are in conditions of \lemref{ind5}. The backward orbit 
(\ref{z-orbit}) either $\eps$-jumps at a moment corresponding to a
return of the orbit (\ref{J-orbit-2}) to $I^{m(\chi)}$, or 
$z''\in D(I^{m(\chi)})$, and 
$z'''\in D(B^{l+1})$. This puts us in the conditions of 
\lemref{ind2} again, and the induction step from level $l$ to $l+1$
is completed.

In the case when $F_l$ is immediately renormalizable, by \lemref{immed}
either the backward orbit (\ref{z-orbit}) $\eps$-jumps at a good moment,
corresponding to a return of the orbit (\ref{J-orbit-2}) to $A^l$, or
$z''\in D(B^{l+1})$, and we can apply the \lemref{ind2}.

 Q.E.D.          end of comment }


\medskip\noindent{\bf  Proof of \lemref{sector}} Let us first show that $J(f_k)\subset D_\theta(S^k)$
with a $\theta=\theta(\bar p)$ (recall that $S^k\ni 0$  is the maximal interval on which $f_k$ is 
unimodal).
  
\comm{By  \corref{contr2}, $\diam J(f_k)\leq C(\bar p)|B^k|$.
Take a $z\in J(f_k)$. Let  the first possibility of \lemref{summary} occur,
and a point $z_{-s}$ $\eps$-jumps. Then $z_{-s}\in D_{\delta}( J_{-s})$ 
with $\delta=\delta(\bar p)>0$, since $\dist (z_{-s}, J_{-s})$ is commensurable
 with $|J_{-s}|$. But then by the Schwarz lemma and \lemref{square root},
$z_{-(n-1)}\in D_{\theta}(J_{-(n-1)})$ with a $\theta=\theta(\bar p)>0$.

If the second possibility of \lemref{summary} happens, we are done immediately.
Q.E.D.
end comment}
By  \corref{contr2}, $\diam J(f_\tau)\leq C(\bar p)|B^\tau|$.
 Take $\zeta''\in J(f_\tau)$.
 Let  $\zeta'=f_\tau(\zeta'')$, $\zeta=f_\tau(\zeta')$,
 and $\zeta=\zeta_0,
\zeta_{-1},\dots,\zeta_{-n}=\zeta',\dots,\zeta_{-2n}=\zeta''$
 be the
 corresponding backward orbit.

Let the first possibility of \lemref{summary} occur and $\zeta_{-s}$ 
$\eps$-jumps at a good moment for $s\leq n-1$.  
Then $\zeta_{-s}\in D_{\delta}( J_{-s})$ 
with $\delta=\delta(\bar p)>0$, since $\dist (\zeta_{-s}, J_{-s})$ is
 commensurable. But then by the Schwarz lemma and \lemref{square root},
and $\zeta''\in Q_\theta (S_\tau)$ with a $\theta=\theta(\bar p)>0$,.

%
Let the second possibility of \lemref{summary} occur.

Let us first consider the case when $f_{\tau -1}$ is not immediately
renormalizable. 
Then $\zeta'\in D(I^{\tau-1,m(\chi-1)})$. By \lemref{ind4},
$\zeta''\in D(I^{\tau-1,m(\chi)})\subset D(S^\tau)$.
Thus $J(f_\tau)\subset  Q_\eps(S_\tau)$, and we are done. 

In the case when $f_{\tau -1}$ is immediately renormalizable $\zeta'\in
 D(B^{\tau -1})$. Consider the interval of monotonicity of $f_{\tau-1}$,
$H=H_{2}(\zeta'')\subset S_{\tau}$. By \lemref{joc-2}, $\zeta''\in
 D_\theta(H)$ with $\theta=\pi/2-O(\eps)$, and the claim follows.

Let us now show how to replace $S^\tau$ by $B^\tau$. By essentially bounded geometry, 
the space $S^\tau\setminus B^\tau$ is commensurable with $|B^\tau|$. 
Also, $S^\tau$ is well inside
$T^\tau=f_\tau S^\tau$. It follows that for any $\delta>0$, there is an $N=N(\bar p, \delta)$ such that
the $N$-fold pull-back of $S^\tau$ by $f_\tau$ is contained in
 $(1+\delta)B^\tau$. 
By the Schwarz lemma and \lemref{square root}, $J(f_\tau)\subset D_{\rho}((1+\delta)B^\tau),$
with a $\rho=\rho(\delta, \bar p)$. 

But for some $\delta>0$ (independent of $\tau$) the map $f_\tau$ 
is linearizable in the $\delta|B^\tau|$-neighborhood 
of the fixed point $\beta_\tau$. In the corresponding  local chart the Julia set $J(f_\tau)$ is invariant with respect to 
$f'_\tau(\beta_\tau)$-dilation.
 Hence further pull-backs will keep it within  a definite sector.  \QED

\end{document}